\documentclass[12pt]{article}
\usepackage{amsfonts,amssymb}
\usepackage{amsmath}
\usepackage{amsthm}

\title{Some open problems}
\author{A.~A.~Agrachev\thanks{SISSA, Trieste \& Steklov Math. Inst.,
Moscow}}
\date{}
\begin{document}
\maketitle

\begin{abstract} We discuss some challenging open problems in the geometric control theory and sub-Riemannian geometry.
\end{abstract}

It is getting harder to prove theorems and easier to force other people to prove them when you are sixty.
Some colleagues asked me to describe interesting open problems in geometric control and sub-Riemannian geometry. Here I list few really challenging problems; some of them are open for a long time and were publicly or privately stated by well-known experts: \linebreak J.-M.~Coron, I.~Kupka, R.~Montgomery, B.~Shapiro, H.~Sussmann, and others.

\bigskip {\bf I. Singularities of time-optimal trajectories.}

\medskip Let $f,g$ be a pair of smooth (i.\,e. $C^\infty$) vector fields on a $n$-dimensional manifold $M$. We study time-optimal trajectories for the system
$$
\dot q=f(q)+ug(q),\quad |u|\le 1,\ q\in M,
$$
with fixed endpoint. Admissible controls are just measurable functions and admissible trajectories are Lipschitz curves in $M$. We can expect more regularity from time-optimal trajectories imposing reasonable conditions on the pair of vector fields.

\smallskip {\sl A. $(f,g)$ is a generic pair of vector fields.} Optimal trajectories cannot be all smooth; are they
piecewise smooth? This is true for $n=2$. More precisely, if $\dim M=2$, then any point of $M$ has a neighborhood such that any contained in the neighborhood time-optimal trajectory is piecewise smooth with atmost 1 switching point (see \cite{Su1,BoPi}). According to the control theory terminology, a switching point of an admissible trajectory is a point where the trajectory is not smooth.

The question is open for $n=3$. What is known? Let $\mathrm{sw(q)}$ be minimal among numbers $k$ such that any contained
in a sufficiently small neighborhood of $q\in M$ time-optimal trajectory has no more than $k$ switching points. We set $\mathrm{sw}(q)=\infty$ if any neighborhood of $q$ contains a time-optimal trajectory with an infinite number of switching points. It is known that $\mathrm{sw}(q)=2$ for any $q$ out of a 2-dimensional Whitney stratified subset of the 3-dimensional manifold $M$ (see \cite{Su2,Sc}) and $\mathrm{sw}\le 4$ for any $q$ out of a 1-dimensional Whitney stratified subset of $M$ (see \cite{AgSi}). Some further results in this direction can be found in \cite{Si}. We do not know if $\mathrm{sw}(q)<\infty$ for any $q\in M$. We also do not know if a weaker property, the finiteness of the number of switching points for any
individual time-optimal trajectory is valid.

Higher dimensions. There is a common opinion that starting from some (not very big) dimension, time-optimal trajectories with accumulating switching points cannot be eliminated by a $C^\infty$-small perturbation of the
system and thus survive any genericity conditions. However, to my knowledge, this opinion was never supported
by a proof. There are very interesting examples of extremals with accumulating switching points whose structure
survives small perturbations (see \cite{Ku,ZeBo}) but nobody knows if these extremals are optimal.

\smallskip {\sl B. $f, g$ are real analytic vector fields.} Let $M$ be a real analytic manifold and $f,g$ analytic vector fields, not necessary generic. Here we cannot expect any regularity of an arbitrary time-optimal trajectory. Indeed, it is possible, even for linear systems, that all admissible trajectories are time-optimal. We can
however expect that among all time-optimal trajectories connecting the same endpoints there is at least one not so bad.

If $n=2$, then any two points connected by a time-optimal trajectory can be connected by a time-optimal trajectory
with a finite number of switching points (see \cite{Su4,Su5}). This is not true for $n\ge 3$. Indeed, classical
Fuller example with accumulating switching points \cite{Fu} can be easily reformulated as a 3-dimensional time-optimal problem. Main open question here is as follows: Given two points connected by a time-optimal trajectory, can we connect them by a time-optimal trajectory with no more than a countable number of switching points?

What is known? The points can be connected by a time-optimal trajectory whose set of switching points is nowhere
dense \cite{Su3}, is not a Cantor set (can be derived from \cite{Ag1}), and satisfies some additional restrictions \cite{Si}. We do not know if we can avoid a positive measure set of switching points.

All mentioned open questions are not easy to answer. In my opinion, the most interesting is one on generic 3-dimensional systems.

\bigskip {\bf II. Cutting the corners in sub-Riemannian spaces.}

\medskip Unlikely the just discussed problems, optimal paths in sub-Riemannian geometry are usually smooth.
However we do not know if they are always smooth. A natural open question here is as follows.
Let $\gamma_i:[0,1]\to M$ be two smooth admissible paths of a sub-Riemannian structure on the manifold $M,\
\gamma_0(0)=\gamma_1(0)=q_0,\ \dot\gamma_0(0)\wedge\dot\gamma_1(0)\ne 0$. Does there exist an admissible path connecting
$\gamma_0(1)$ with $\gamma_1(1)$ that is strictly shorter than the concatenation of the curves $\gamma_0$ and $\gamma_1$?

Admissible paths are integral curves of a bracket generating vector distribution $\Delta\subset TM$. It is easy to
show that positive answer to the question for rank 2 distributions implies positive answer in the general case.
Let $\Delta=span\{f_0,f_1\}$, where $f_0,f_1$ are smooth vector fields on $M,\ \dim M=n$. We set
$$
n_k(q)=span\left\{[f_{i_1},[\cdots,f_{i_j}]\cdots](q): i_j\in\{0,1\},\ j\le k\right\},
$$
$m=\min\{k: n_k(q_0)=n\}.$
If $m\le 4$, then the answer to our question is positive: it is proved in \cite{LeMo}. Moreover, an example
studied in \cite{Mo} supports the conjecture that the answer is perhaps positive for $m=5,\ n\le 4$ as well.
Any improvement of the estimates for $m$ and $n$ would be very interesting. We still know very little about
sub-Riemannian structures with big $m$ and it may happen that the answer is negative for some $m$ and $n$.

\bigskip{\bf III. ``Morse--Sard theorem" for the endpoint maps.}

\medskip We continue to consider admissible paths of a sub-Riemannian structure on $M$. Given $q_0\in M$,
the space of starting at $q_0$ admissible paths equipped with the $H^1$-topology forms a smooth Hilbert
manifold. The endpoint map is a smooth map from this Hilbert manifold to $M$; it sends a path
$\gamma:[0,1]\to M$ to the point $\gamma(1)$. Critical points of the endpoint map are called {\it singular
curves} of the distribution.

The Morse--Sard theorem for a smooth map defined on a finite dimensional manifold states that the set of
critical values of the map has zero measure. It is not true in the infinite dimensional case: there are
smooth surjective maps without regular points from any infinite dimensional Banach space to $\mathbb R^2$ (see \cite{Ba}).

The endpoint maps have plenty of regular points but we do not know if they always have regular values.
This is an interesting open question. We can reformulate the question as follows: is it possible that
starting from $q_0$ singular curves fill the whole manifold $M$?

Optimal (i.\,e. length minimizing) singular curves are better controlled; we know that starting from $q_0$
optimal singular curves fill a nowhere dense subset of $M$ (see \cite{Ag2}). An important open question:
can they fill a positive measure subset of $M$?

\bigskip{\bf IV. Unfolding the sub-Riemannian distance.}

\medskip The problem concerns singularities of the distance function for generic sub-Riemannian structures.
Let $q_0\in M$ and $S_{q_0}: M\to[0,+\infty)$ be the sub-Riemannian distance from the point $q_0$.
Sufficiently small balls $S^{-1}\left([0,\varepsilon]\right)$ are compact. Let $q\in M$ be a point from
such a ball. Then $q$ is connected with $q_0$ by an optimal path. If this optimal path is not a singular curve and any point from a neighborhood of $q$ is connected with $q_0$ by a unique optimal path, then $S_{q_0}$
is smooth at $q$.

The points connected with $q_0$ by more than one optimal path form the {\it cut locus}. The function $S_{q_0}$
is not smooth in the points of the cut locus and it is not smooth at the points connected with $q_0$ by
optimal singular curves but these two types of singularities are very different.

If all connecting $q_0$ and $q$ optimal paths are not singular curves, then the singularity of $S_{q_0}$
at $q$ is similar to singularities of Riemannian distances and, more generally, to singularities of the
optimal costs of regular variational problems. The function $S_{q_0}$ is semiconcave \cite{CaRi} and typical
singularities in low dimensions are well-described by the theory of Lagrangian and Legendrian singularities
\cite[Ch.\,3]{ArVaGu} developed by V.~Arnold and his school.

 On the other hand, if $q_0$ is connected with $q$ by an optimal singular curve, then $S_{q_0}$ is not even
 locally Lipschitz at $q_0$ (see \cite[Ch.\,10]{AgBaBo}); moreover, classical singularities theory does not work
 and the structure of typical singularities is totally unknown.
 There are few studied models \cite{AgBoChKu,Sa1} but they are too symmetric to be typical and the structure of
 their singularities is easily destroyed by small perturbations.

 Let us consider, in particular, the {\it Martinet distribution} that is a rank 2 distribution in $\mathbb R^3$
 in a neighborhood of a point $q_0$ such that $n_2(q_0)=2,\ n_3(q_0)=3$ (see {\bf II.} for the definition
 of $n_i(q_0)$). The points $q$ in a neighborhood of $q_0$ where $n_2(q)=2$ form a smooth
 2-dimensional submanifold $N\subset M$, the {\it Martinet surface}. Moreover, the distribution $\Delta$ is
 transversal to $N$ and $\Delta_q\cap T_qN,\ q\in N$ is a line distribution on $N$. Integral curves of this
 line distribution are singular curves whose small segments are optimal. There are no other singular curves for
 such a distribution.

 {\sl Example.} Let $f_1=\frac\partial{\partial x_1},\
 f_2=\frac\partial{\partial x_2}+x_1^2\frac\partial{\partial x_3}$, then $\Delta=span\{f_1,f_2\}$ is a
 Martinet distribution and the Martinet surface is a coordinate plane defined by the equation $x_1=0$.
 The fields $f_1,f_2$ form an orthonormal frame of the so called `flat' sub-Riemannian metric on the Martinet
 distribution. Let $q_0=0$, singularities of $S_0$ are well-known (see \cite{AgBoChKu}). The cut locus has the
 form: $\{x\in\mathbb R^3: x_1=0,\ x_2\ne 0\}$, the Martinet surface with the removed singular curve through
 $q_0$. The singular locus of a sphere $S_{q_0}^{-1}(\varepsilon)$ is a simple closed curve and its complement
 (the smooth part of the sphere) is diffeomorphic to the disjoint union of two discs.

 The `flat' metric is rather symmetric, in particular, it respects the orthogonal reflection of $\mathbb R^3$
 with respect to the Martinet plane. Simple topological arguments show that for generic metric with a broken
 symmetry, the smooth part of a sphere is connected and is not contractible. The singular locus of the sphere
 should be cut at the points where the sphere intersects the optimal singular curve but the shape of the sphere near these points is unknown.

 An important open question is to find a $C^1$-classification of the germs of spheres at the points of optimal
 singular curves for generic metrics. Here we say that two germs are $C^1$-equivalent if one can be
 transformed into another by a germ of $C^1$-diffeomorphism of $\mathbb R^3$.

 The next step is the Engel distribution, i.\,e. a rank 2 distribution in $\mathbb R^4$ such that $n_3(q)=3,\
 n_4(q)=4$. There is exactly one singular curve through any point and small segments of singular curves are
 optimal. We repeat our question for this case; the spheres are now 3-dimensional hypersurfaces in $\mathbb R^4$.

Generic germ of a rank 2 distribution in $\mathbb R^n$ possesses a $(n-4)$-dimensional family of singular curves
through $q_0\in\mathbb R^n$ (see, for instance, \cite{Mon}). Take a generic curve from this family; its small
segments are optimal. Take a point $q$ where the selected singular curve intersects the sphere
$S_{q_0}^{-1}(\varepsilon)$. The points of singular curves from the family in a small neighborhood of $q$ in our sphere form a smooth $(n-4)$-dimensional submanifold $\Sigma\subset\mathbb R^n$. The intersection of the sphere
with a transversal to $\Sigma$ smooth 4-dimensional submanifold should have a shape similar to the germ of the
sphere in the Engel case. A neighborhood of $q$ in the sphere is fibered by such intersections. Hence the solution
of the problem in the 4-dimensional Engel case is a very important step in the unfolding of the sphere for any
$n\ge 4$.

The desired classification seems to be complicated. There is a 2-dimensional modification of the problem that,
in my opinion, already contains the main difficulty. Ones resolved, it will reduce the study of higher dimensional
problems to the conventional singularities theory techniques. Consider the germ at $q_0\in\mathbb R^2$ of a pair
of smooth vector fields $f_0,f_1$ such that $f_0(q_0)\wedge f_1(q_0)=0,\ n_2(q_0)=1,\ n_3(q_0)=2$. The almost
Riemannian distance $S_{q_0}(q)$ is the optimal time to get $q$ from $q_0$ by an admissible trajectory of the system
$$
\dot q=u_0f_0(q)+u_1f_1(q), \quad u_1^2+u_2^2=1.
$$
The question is to find a $C^1$-classification of the germs of distance functions $S_{q_0}$ for generic pairs of
vector fields $f_0,f_1$ among the pairs that satisfy conditions $n_2(q_0)=1,\ n_3(q_0)=2$. See \cite{Bon}
for some partial results.

\bigskip{\bf V. Symmetries of vector distributions.}

\medskip A symmetry of a distribution $\Delta\subset TM$ is a diffeomorphism $\Phi:M\to M$ such that
$\Phi_*\Delta=\Delta$. The differential geometry appeals to search most symmetric objects in the class, those with
a maximal symmetry group. The singularities theory, on the contrary, encourages the study of less symmetric generic objects. Both paradigms have their reasons and complement each other. Anyway, a fundamental problem is to characterize objects whose symmetry groups are finite-dimensional Lie groups.

Our objects are vector distributions. Any symmetry transfers singular curves of the distribuiton in singular
curves and these curves often play a key role in the calculation of symmetry groups (see \cite{KrZe,DuZe}).
We say that a distribution is singular transitive if any two points of $M$ can be connected by a concatination
of singular curves. A natural open question is as follows: Is it true that singular transitivity of the
distribution implies that its symmetry group is a finite dimensional Lie group?

All known examples support the positive answer to this question. Moreover, the group of symmetries
is infinite dimensional for many popular classes of not singular transitive distribution: codimension 1
distributions, involutive distributions, Goursat--Darboux distributions. We can even expect that any not singular
transitive rank 2 distribution has an infinite dimensional symmetry group. Some results of \cite{DuZe}
seem to be rather close to this statement.
\newpage

{\bf VI. Closed curves with a nondegenerate Frenet frame.}

\medskip Let $\gamma:S^1\to\mathbb R^n$ be a smooth closed curve in $\mathbb R^n$. We say that $\gamma$ is
degenerate at $t\in S^1$ if $\dot\gamma(t)\wedge\cdots\wedge\gamma^{(n)}(t)=0$. Degeneracy points are the points
where velocity or curvature of the curve vanishes if $n=2$, where velocity or curvature or torsion vanishes if $n=3$ e.\,t.\,c. The curve is nondegenerate if it has no degeneracy points. Any nondegenerate curve admits the
orthonormal Frenet frame $E(t)=\left(e_1(t),\ldots,e_n(t)\right),\ t\in S^1$, that
is a smooth closed curve in the orthogonal group $\mathrm O(n)$.

Now let $n=3$ and $\gamma$ be a plane convex curve, $\gamma(t)\in\mathbb R^2\subset\mathbb R^3,\
\forall t\in S^1$. Then any small perturbation of $\gamma$ as a spatial curve is degenerate in some points.
On the other hand, an appropriate small perturbation of a plane convex curve run twice (say, of the curve
$t\mapsto\gamma(2t),\ t\in S^1$) makes it a nondegenerate curve in $\mathbb R^3$. Everyone can get evidence of
that playing with a cord on the desk. This is also a mathematical fact proved in \cite{Fe,Mi}.

Frenet frame of the plane convex curve treated as a spatial curve is a one-parametric subgroup
$\mathrm{SO}(2)\subset\mathrm O(3)$, a shortest closed geodesic in $\mathrm O(3)$ equipped with a standard
bi-invariant metric. It is proved in \cite{Mi} that the length of the Frenet frame of any regular curve in
$\mathbb R^3$ is greater than the double length of $\mathrm{SO}(2)$.

Come back to an arbitrary $n$. Let $\mu(n)$ be minimal $m$ such that run $m$ times convex plane curves have
regular small perturbations in $\mathbb R^n$. We know that $\mu(2)=1,\ \mu(3)=2$. An important open problem
is to find $\mu(n)$ for $n>3$ and to check if the length of the Frenet frame of any regular curve in $\mathbb R^n$
is greater than the length of $\mathrm{SO}(2)\subset\mathrm O(n)$ multiplied by $\mu(n)$.

Let me explain why this problem is a challenge for the optimal control theory and why its study may bring important new tools to the theory. The Frenet structural equations for a regular curve $\gamma$ in $\mathbb R^n$ have a form:
$$
\dot\gamma=e_1,\quad \dot e_i=u_k(t)e_{i+1}-u_{i-1}(t)e_{i-1},\ i=1,\ldots,n-1, \eqno (1)
$$
where $u_0=u_n=0,\ u_i(t)>0,\ i=1,\ldots,n-1,\ t\in S^1$.

In other words, regular curves together with there Frenet frames are periodic admissible trajectories of the control system (1) with positive control parameters $u_1,\ldots,u_n$. The length of the Frenet frame on the
segment $[0,t_1]$ is $\int\limits_0^{t_1}\left(u_1^2(t)+\cdots+u_{n-1}^2(t)\right)^{\frac 12}dt$.
We are looking for a periodic trajectory with shortest Frenet frame.

 A shortest frame is unlikely to exists since control parameters belong to an open cone. It is reasonable to expect that minimizing sequences converge to a solution of (1) with
$u_2(t)\equiv\cdots\equiv u_{n-1}(t)\equiv 0$, while $u_1(t)$ stays positive to guarantee the periodicity of $\gamma$. In other words, the infimum is most likely realized by a plane convex curve run several times.
Obviously, the length of the Frenet frame does not depend on the shape of the convex curve.

So we have to take the $m$ times run circle: $u_1(t)=1,\ u_2(t)=\cdots=u_{n-1}(t)=0,\ 0\le t\le2\pi m$,
and try to find small positive perturbations of control parameters in such a way that the perturbed curve
stays periodic. Then $\mu(n)$ is minimal among $m$ for which such a perturbation does exist. Unfortunately,
we cannot use typical in geometric control sophisticated two-side variations that produce iterated Lie brackets:
only one-side variations are available. I think, it is a very good model to understand high order effects
of time-distributed one-side variations.

The study of the 3-dimensional case by Milnor in \cite{Mi} was not variational; it was a nice application of
the integral geometry. However, the integral geometry method is less efficient in higher dimensions (see
\cite{NoYa} for some partial results).

\bigskip{\bf VII. Controllability of the Navier--Stokes equations controlled by a localized degenerate forcing.}

\medskip We consider the Navier--Stokes equation of the incompressible fluid:
$$
\frac{\partial u}{\partial t}+(u,\nabla)u-\nu\Delta u+\nabla p=\eta(t,x),\quad \mathrm{div}u=0, \eqno (2)
$$
with periodic boundary conditions: $x\in\mathbb T^d/2\pi\mathbb Z^d,\ d=2,3.$
Here $u(t,x)\in\mathbb R^d$ is the velocity of the fluid at the point $x$ and moment $t;\ \nu$ is a positive
constant (viscosity), $p$ is the pressure and $\eta$ external force.

We treat (2) as an evolution equation in the space of divergence free vector fields on the torus $\mathbb T^d$
controlled by the force. In other words, $u(t,\cdot)$ is the state of our infinite dimensional control
system and $\eta$ is a control. These notations are against the control theory tradition where $u$ is always
control but we do not want to violate absolutely standard notations of the mathematical fluid dynamics.
By the way, symbol $u$ for the control was introduced by Pontryagin as the first letter of the Russian word
``upravlenie" that means control.

The state space is $V=\left\{u\in H^1(\mathbb T^d,\mathbb R^d): \mathrm{div}u=0\right\}$, control parameters
$\eta(t,\cdot)$ belong to a subspace $E\subset V$. We say that the system is approximately controllable
$\bigl($controllable in finite dimensional projections$\bigr)$ in any time if  for any $u_0,u_1\in V,\ t_1>0$
and any
$\varepsilon>0$ $\bigl($any finite dimensional subspace $F\subset V\bigr)$ there exists a bounded control
$\eta,\ \eta(t,\cdot)\in E,\ 0\le t\le t_1$, and a solution $u$ of (2) such that $u(0,\cdot)=u_0,\
\|u(t_1,\cdot)-u_1\|_{L_2}<\varepsilon\quad \bigl(P_F(u(t_1,\cdot)-u_1)=0$, where $P_F$ is the $L_2$-orthogonal
projector on $F\bigr)$.

Of course, controllability properties depend on the choice of the space of control parameters $E$. It is known
that the systems is controllable in both senses by a {\it localized forcing} when
$E=\{u\in V: \mathrm{supp}\,u\subset\bar{\mathcal D}\}$ and $\mathcal D$ is an arbitrary open subset of
$\mathbb T^d$. Moreover, such $E$ provides a much stronger exact controllability (see
\cite{Co,CoFu,FuIm,Im}).

On the other hand, the system is approximately controllable and controllable in finite-dimensional projections
by a {\it degenerate forcing} (or forcing with a localized spectrum) when $E$ is a finite dimensional
space of low frequency trigonometric polynomials (see \cite{AgSa1,AgSa2,Sh1,Sh2}). This kind of controllability
illustrates a mechanism of the energy propagation from low to higher frequencies that is a necessary step
in the long way towards a reliable mathematical model for the well-developed turbulence.

It is important that the control parameters space $E$ does not depend on the viscosity $\nu$. Moreover, if
$d=2$, then the described controllability properties are valid also for the Euler equation (i.\, e. for $\nu=0$);
the Cauchy problem for the Euler equation is well-posed in this case.

Now an important open question: is the system approximately controllable and (or) controllable in the finite
dimensional projections by a {\it localized degenerate forcing} when $E$ is a finite dimensional subspace of
the space $\{u\in V: \mathrm{supp}\,u\subset\bar{\mathcal D}\}$? The question is about existence and effective
construction of such a space $E$ that does not depend on the viscosity $\nu$.

The independence on $\nu$ is important for eventual applications to the well-developed turbulence that concerns
the case of very small $\nu$ (or very big Reynold number). Of course, similar problems for other boundary
conditions and other functional spaces are also very interesting.

\bigskip We have arrived to a sacral number of seven problems and can relax a little bit. To conclude,
I would like
to discuss one more problem; it is less precise than already stated questions but, to my taste, is nice and fascinating. The problem concerns contact 3-dimensional manifolds and is inspired by the Ricci flow story.

\newpage
{\bf $\infty$. Diffusion along the Reeb field.}

\medskip I recall that a contact structure on a 3-dimensional manifold $M$ is a rank 2 distribution $\Delta\subset M$ such that $n_2(q)=3,\ \forall q\in M$. According to a classical Martinet theorem, any orientable 3-dimensional
manifold admits a contact structure. I am going to introduce some dynamics on the space of sub-Riemannian metrics
on a fixed compact contact manifold $(M,\Delta)$.

First, to any sub-Riemannian metric on $(M,\Delta)$ we associate a transversal to $\Delta$ Reeb vector field $e$ on $M$. In what follows, we assume that $\Delta$ is oriented; otherwise $e(q)$ is defined up-to a sign but
further considerations are easily extended to this case. Let $\omega$ be a nonvanishing differential 1-form
on $M$ that annihilates $\Delta$. The condition $n_2(q)=3$ is equivalent to the inequality
$\omega_q\wedge d_q\omega\ne 0$. The form $\omega$ is defined up-to the multiplication by a nonvanishing function; the sign of the 3-form $\omega_q\wedge d_q\omega$ does not depend on the choice
of $\omega$ and defines an orientation on $M$. We have: $d_q\omega\bigr|_{\Delta_q}\ne 0$; moreover, $d_q(a\omega)\bigr|_{\Delta_q}=a(q)d_q\omega\bigr|_{\Delta_q}$ for any smooth function $a$ of $M$.

Given a
sub-Riemannian metric on $\Delta$, there exists a unique annihilating $\Delta$ form $\omega$ such that the
2-form $d_q\omega\bigr|_{\Delta_q}$ coincides with the area form on $\Delta_q$ defined by the inner product
and the orientation. The kernel of $d_q\omega$ is a 1-dimensional subspace of $T_qM$ transversal to $\Delta_q$,
and $e(q)$ is an element of this kernel normalized by the condition $\langle\omega_q,e(q)\rangle=1$.

In other words, the Reeb vector field is defined by the conditions: \linebreak $i_e\omega=1,\ i_ed\omega=0$. Hence
$L_e\omega=0$, where $L_e$ is the Lie derivative along $e$, and the generated by $e$ flow on $M$ preserves $\omega$. In general, this flow does not preserves the sub-Riemannian metric.
We may try to classify contact structures by selecting best possible sub-Riemannian metrics on them.

Assume that there exists a metric preserved by the flow generated by the Reeb vector field. Take a standard extension of the sub-Riemannian metric to a Riemannian metric on $M$: simply say that $e$ is orthogonal to
$\Delta$ and has length 1. The generated by $e$ flow preserves this Riemannian metric as well. So our
compact Riemannian space admits a one-parametric group of isometries without equilibria. Hence $M$ is a
Seifert bundle. Do not care if you do not remember what is Seifert bundle: it is sufficient to know
that they are classified as well as invariant contact structures on them.

The invariant with respect to the Reeb field sub-Riemannian metric gives a lot of information about the
manifold. Let $q\in M$; our sub-Riemannian metric induces a structure of Riemannian surface on a neighborhood
of $q$ factorized by the trajectories of the local flow generated by the restriction of $e$ to the neighborhood.
Let $\kappa(q)$ be the Gaussian curvature of this Riemannian surface at the point $q$; then $\kappa$ is a
well-defined smooth function on $M$, a differential invariant of the sub-Riemannian metric. Moreover, $\kappa$
is a first integral of the flow generated by the Reeb field $e$. If $\kappa=0$, then universal covering of the
sub-Riemannian manifold is isometric to the Heisenberg group endowed with the standard left-invariant metric.
If $\kappa$ is a negative (positive) constant, then universal covering of the sub-Riemannian manifold is isometric to the universal covering of the group $\mathrm{SL}(2)$ (group $\mathrm{SU}(2)$) equipped with a left-invariant
sub-Riemannian metric induced by the Killing form.

Assume that function $\kappa$ is not a constant and $c\in\mathbb R$ is it regular value. Then $\kappa^{-1}(c)$
is a compact 2-dimensional submanifold of $M$; we treat it as a 2-dimensional Riemannian submanifold of the
Riemannian manifold $M$ equipped with the standard extension of the sub-Riemannian structure. It is easy to
see that $\kappa^{-1}(c)$ is isometric to a flat torus. Indeed, $T(\kappa^{-1}(c))$ contains the field
$e|_{\kappa^{-1}(c)}$ and is transversal to $\Delta$; the field $e|_{\kappa^{-1}(c)}$ and the
unit length field from the line distribution $T(\kappa^{-1}(c))\cap\Delta$ commute and form an orthonormal
frame.

So preserved by the Reeb field sub-Riemannian metrics have plenty of nice properties. Unfortunately, not
any compact contact manifold admits such a metric because not any compact 3-dimensional manifold admits
a structure of Seifert bundle. I am going to discuss a natural procedure that may lead to a generalized
version of such a metric with reasonable singularities.

It is more convenient to work in the cotangent bundle than in the tangent one. A sub-Riemannian metric
is an inner product on $\Delta\subset TM$; let us consider the dual inner product on
$\Delta^*=T^*M/\Delta^\perp$, where $\Delta^\perp$ is the annihilator of $\Delta$. This is a family of
positive definite quadratic forms on $\Delta^*=T_q^*M/\Delta_q^\perp,\ q\in M$, or, in other words, a family
of nonnegative quadratic forms $h_q$ on $T^*_qM$ such that $\ker h_q=\Delta_q^\perp$. The function
$$
h:T^*M\to\mathbb R,\ \mathrm{where}\ h(\xi)=h_q(\xi),\quad \forall\xi\in T^*_qM,\ q\in M,
$$
is the {\it Hamiltonian of the sub-Riemannian metric}. Hamiltonian vector field on $T^*M$ associated to $h$
generates the sub-Riemannian geodesic flow.

The Hamiltonian $h$ determines both the vector distribution and the inner product. We denote by
$u_h:T^*M\to\mathbb R$ the Hamiltonian lift of the Reeb field $e$,
$$
u_h(\xi)=\langle\xi,e(q)\rangle,\quad \forall\xi\in T^*_qM,\ q\in M,
$$
and by $\mathcal U^t_h:T^*M\to T^*M,\ t\in\mathbb R,$ the Hamiltonian flow generated by the Hamiltonian
field associated to $u_h$. The flow $\mathcal U^t_h$ is a lift to the cotangent bundle of the flow on
$M$ generated by $e$. Let $P_t:M\to M$ be such a flow, $\frac{\partial P_t(q)}{\partial t}=e\circ P_t(q),\
P_0(q)=q,\ q\in M$; then $\mathcal U^t_h=P^*_{-t}$. The flow $P_t$ preserves the sub-Riemannian metric if and only
if the flow $\mathcal U^t_h$ preserves $h$; in other words, if and only if $\{u_h,h\}=0$, where
$\{\cdot,\cdot\}$ is the Poisson bracket. Note that $h\bigr|_{T^*_qM}$ is a quadratic form and
$u_h\bigr|_{T^*_qM}$ is a linear form, $\forall q\in M$; hence $\{u_h,h\}\bigr|_{T^*_qM}$ is a quadratic form.

Recall that the flow $\mathcal U^t_h=P^*_{-t}$ preserves the 1-form $\omega$, and $\omega$ is a nonvanishing
section of the line distribution $\Delta^\perp$. Hence $\Delta^\perp$ is contained in the kernel of the
quadratic forms $h\circ\mathcal U^t_h\bigr|_{T^*_qM}$ and
$\underbrace{\{u_h,\{\cdots\{u_h}_i,h\}\cdots\}\bigr|_{T^*_qM}=
\frac{d^i}{dt^i}\bigr|_{t=0}\left(h\circ\mathcal U^t_h\right)\bigr|_{T^*_qM}$.

We are now ready to introduce the promised dynamics on the space of sub-Riemannian metrics on $\Delta$, where
metrics are represented by their Hamiltonians. Let $\varepsilon$ be a positive smooth function on $M$. A discrete time dynamical system transforms a Hamiltonian $h_n$ into the Hamiltonian
$$
h_{n+1}=\frac 1{2\varepsilon}\int\limits_{-\varepsilon}^\varepsilon h_n\circ\mathcal U^t_{h_n}\,dt,\quad n=0,1,2,\ldots,
$$
a partial average of $h_n$ with respect to the flow $\mathcal U^t_h$.

The Hamiltonian $h_{n+1}$ is
equal to $h_n$ if and only if $\{u_{h_n},h_n\}=0$.
Indeed, let $\langle\cdot,\cdot\rangle_q$ be an inner
product in $\Delta^*_q$ and $H^t_q:\Delta^*_q\to\Delta^*_q$ the symmetric operator associated to the quadratic
form $h_n\circ\mathcal U^t_{h_n}\bigr|_{\Delta^*_q}$ by this inner product: $h_n\circ\mathcal U^t_{h_n}(\cdot)=
\langle H^t_n\cdot,\cdot\rangle_q$. Recall that the flow $\mathcal U^t_{h_n}$ is generated by the Reeb field of
$h_n$, hence the area form on $\Delta^*_q$ defined by $h_n\circ\mathcal U^t_{h_n}\bigr|_{\Delta^*_q}$ does not
depend on $t$; in other words, $\det H^t_n=const$. The equation $\det H=const$ defines a strongly convex
hyperboloid in the 3-dimensional cone of positive definite symmetric operators on the plane, and $H^t_n$ is a
curve in such a hyperboloid; hence $\frac 1{2\varepsilon}\int\limits_{-\varepsilon}^\varepsilon H^t_n\,dt=H^0_n$ if and only if $H^t_\varepsilon\equiv H^0_n$, i.\,e. $h_n\circ\mathcal U^t_{h_n}\equiv h_n$.

If the sequence $h_n$ converges, then its limit is the Hamiltonian of a sub-Riemannian metric on $\Delta$
preserved by the Reeb field. Otherwise we may modify the sequence and take scaled averages:
$$
h_{n+1}=c_n\int\limits_{-\varepsilon_n}^{\varepsilon_n} h_n\circ\mathcal U^t_{h_n}\,dt.
$$
There is a good chance to arrive to a nonzero limiting Hamiltonian $h_\infty$ by a clever choice of the sequences
of positive functions $\varepsilon_n,c_n$. Then
$h_\infty\bigr|_{T^*_qM}$ is a nonnegative quadratic form and $\Delta^\perp_q\subset\ker h_\infty\bigr|_{T^*_qM}$
for any $q\in M$. It may happen however that $rank\left(h_\infty\bigr|_{T^*_qM}\right)<2$
for some $q\in M$ and $h$ is not the Hamiltonian of a contact sub-Riemannian metric; we can treat it as a generalized version of such a metric.

A continuous time analogue of the introduced dynamics is a ``heat along the Reeb field" equation
$$
\frac{\partial h}{\partial t}=c\{u_h,\{u_h,h\}\}
$$
in the space of sub-Riemannian metrics on the given contact distribution. It is easy to show that the
equality $\{u_h,\{u_h,h\}\}=0$ implies $\{u_h,h\}=0$ and stationary solutions of this equation are
exactly the metrics preserved by the Reeb fields.
I conclude with an explicit expression for this nice and mysterious evolution equation in the appropriate frame.

All contact distributions are locally equivalent according to the Darboux theorem. Let $f_1,f_2$ be a basis
of the contact distribution $\Delta$ such that $f_1,f_2$ generate a Heisenberg Lie algebra:
$[f_1,[f_1,f_2]]=[f_2,[f_2,f_1]]=0$. We set $v_i(\xi)=\langle\xi,f_i(q)\rangle,\ \xi\in T^*_qM,\ q\in M$, the
Hamiltonian lift of the field $f_i,\ i=1,2$; then
$$
\{v_1,\{v_1,v_2\}\}=\{v_2,\{v_2,v_1\}\}=0. \eqno (3)
$$
Hamiltonian of any sub-Riemannian metric on $\Delta$ has a form:
$$
h=a_{11}v_1^2+2a_{12}v_1v_2+a_{22}v_2^2,  \eqno (4)
$$
where $a_{ij}$ are smooth functions on the domain in $M$ where $f_1,f_2$ form a basis of $\Delta$, and the
quadratic form defined by the matrix $A(q)=\left(\begin{smallmatrix} a_{11}(q)& a_{12}(q)\\
a_{12}(q) & a_{22}(q) \end{smallmatrix}\right)$ is positive definite for any $q$ from this domain.
Let $\delta=\det A$; this is a function on $M$ and we treat it as a constant on the fibers function on $T^*M$.
A key for us function $u_h$ depends only on $\delta$ and has  a form:
$$
-u_h=\delta\{v_1,v_2\}+v_1\{v_2,\delta\}+v_2\{\delta,v_1\}.  \eqno (5)
$$
The relations (3)--(5) give an explicit expression for the equation
$\frac{\partial h}{\partial t}=c\{u_h,\{u_h,h\}\}$ as a system of third order partial differential equations
for the functions $a_{ij}$.

\end{document}